\documentclass[preprint,12pt]{elsarticle}

\usepackage{amssymb}

\usepackage{graphicx}
\usepackage{latexsym}
\usepackage{amsfonts}
\journal{Journal of Mathematical Analysis and Applications}
\newtheorem{assumption}{Assumption}[section]
\newtheorem{lemma}{Lemma}[section]
\newtheorem{theorem}{Theorem}[section]
\newtheorem{problem}{Problem}[section]
\newtheorem{remark}{Remark}[section]
\newtheorem{definition}{Definition}[section]
\newcommand{\Frac}[2]{\displaystyle\frac{#1}{#2}}
\newcommand{\delfrac}[2]{\displaystyle\frac{\partial #1}{\partial #2}}
\newcommand{\iprod}{\, \mbox{\raise.3ex\hbox{\tiny $\bullet$}}\, }
\newcommand{\diffrac}[2]{\displaystyle\frac{d #1}{d #2}}
\newcommand{\QED}{$\Box$} 
\begin{document}

\begin{frontmatter}



\title{Stabilities of Shape Identification Inverse Problems in 
a Bayesian Framework}


\author{Hajime Kawakami}

\ead{kawakami@math.akita-u.ac.jp}

\address{Mathematical Science Course,
Akita University, Akita 010-8502, Japan}

\begin{abstract}
A general shape identification inverse problem is studied 
in a Bayesian framework. 
This problem requires the determination of the unknown shape of a 
domain in the Euclidean space 
from finite-dimensional observation data with some 
Gaussian random noise. 
Then, the stability of posterior is studied for 
observation data.  
For each point of the space, 
the conditional probability that the point is 
included in the unknown domain given the observation data is considered.
The stability is also studied 
for this probability distribution.

As a model problem for our inverse problem,
a heat inverse problem is considered. 
This problem requires the determination of the unknown shape of cavities
in a heat conductor from temperature data of some portion 
of the surface of the heat conductor.  
To apply the above stability results to this model problem, 
one needs  
the measurability and some boundedness of the forward operator.
These properties are shown. 
\end{abstract}


\end{frontmatter}
Declarations of interest: none




\section{Introduction}

In this paper, we study a general shape identification inverse problem.
This problem requires the determination of unknown shape of a 
domain $D$ in $\mathbb{R}^d$ 
from finite-dimensional observation data 
$y \in \mathbb{R}^m,$
that is,
\[
y = F(D) + \eta.
\]
Denote by $\mathcal{D}$ 
the set of all possible domains. 
Then, $F$ is the ``forward'' operator taking one instance 
of input $D \in \mathcal{D}$ into a set
of observations (an $m$-dimensional vector), 
and $\eta$ is an $m$-dimensional Gaussian random noise variable.
For this problem, we use the Bayesian formulation 
given by \cite{DS} and \cite{S1}. 
Although the set of possible domains is not a function space of the same
type as in their theory,
we can consider our problem as based on this previous work   
as mentioned in subsection \ref{subsect.problems.framework}.   
This gives a probability distribution on $\mathcal{D}$  
corresponding to the statistical properties of 
the observational noise.   
 
We can consider many inverse problems as model problems of 
our general shape identification inverse problem.
In this study, we consider the following model problem.
Let $\Omega$ be a bounded domain 
in $\mathbb{R}^d$  
and let $C$ be the union of several (one or more) disjoint connected  
domains in $\Omega$ such that $\overline{C} \subset \Omega.$
Then, $D := \Omega \setminus \overline{C}$ 
is a heat conductor with cavities $C,$
and the shape of $\partial C$ is unknown.   
We consider the following initial 
and mixed boundary value problem for the heat equation: 
\begin{equation}
\left\{
\begin{array}{ll}
\displaystyle{\frac{\partial u}{\partial t}}(t,x)
= \displaystyle{\frac{1}{2}} \Delta u(t,x) 
& \mbox{if $(t,x) \in [0,T] \times D,$}
\\
[0.2cm] 
\displaystyle{\frac{\partial u}{\partial \nu}}(t,x) 
=
\psi(t,x)
& \mbox{if $(t,x) \in [0,T] \times \partial \Omega,$}
\\
u(t,x)=0 & \mbox{if $(t,x) \in [0,T] \times \partial C,$}
\\
u(0,x)=0 & \mbox{if $x \in D.$} 
\end{array} 
\right.
\label{eq.sect.intro.1}
\end{equation}
Here, $T$ is a positive real number, 
$\Delta$ is the standard Laplacian,
$\nu$ is the unit normal to $\partial \Omega$ directed into the 
exterior of $\Omega,$ and 
$\psi$ is a given function on $[0,T] \times \partial \Omega.$ 
Then, $F(D)$ is a finite sampling of values of $u$ on a subset $A$  
of $\partial \Omega = \partial D \setminus \partial C.$ 
We can consider $F = G_2 \circ G_1,$ 
where $G_1(D)$ are the values of $u$ on $[0,T] \times A$ and 
$G_2$ is a finitization operator. 
Then, the inverse problem for $G_1(D)$ is the problem 
of finding $D$ from the observation data $G_1(D).$ 
Many researchers have studied such inverse problems 
in cases where the boundary condition is a Dirichlet condition, 
Neumann condition, Robin condition, 
or mixed condition, e.g., 
\cite{BCSV}, \cite{BC1}, \cite{BC2}, \cite{CKY}, \cite{CRV}, 
\cite{HT}, \cite{HNW}, \cite{IK}, 
\cite{K}, \cite{KT}, and \cite{NW}.

In general shape identification inverse problems  
and in our model problem, the regularity (i.e., smoothness) 
of the unknown shape is important. 
We intend to proceed with discussions 
under as general an assumption as possible.
We consider our model problem 
under the assumption that the unknown shape is Lipschitz continuous. 
\cite{KT} showed the uniqueness for $G_1(D)$ of our model problem
for the case where the unknown shape is Lipschitz continuous. 
\cite{BCSV} and \cite{CRV}  
showed such inverse problems have only the logarithmic stability 
for the observational data 
in the case of a Dirichlet problem, and a Robin problem,
for the case where the unknown shape is slightly smoother.
(The logarithmic stability means that the Hausdorff distance between 
domains can be estimated by some negative power of 
the logarithm of the $L^2$ distance between data boundary value functions.)
This logarithmic stability gives the largest possible error of the estimates. 

However, in the Bayesian formulation,
we can obtain the (locally) Lipschitz stability 
of the probability distribution  
of unknown domains 
for our general shape identification inverse problem 
and our model problem, 
even if the observation data are a finite set and 
the inverse problem is underdetermined. 
In addition to the probability distribution on $\mathcal{D},$ 
we consider for each point $x,$ 
the conditional probability that the point $x$ is 
included in the unknown domain given the observation data $y,$
and we call this probability the ``\textit{domain ratio}''.
We also study the stability of the domain ratio in the 
Bayesian framework. 

The Bayesian formulation 
for our inverse problem is 
given by 
\[
P(\mbox{domain}|\mbox{data}) 
\propto P(\mbox{data}|\mbox{domain})
P(\mbox{domain}),
\]
and the posterior, $P(\mbox{domain}|\mbox{data}),$  
is the probability distribution on $\mathcal{D}.$ 
In our framework, 
some properties of $F(D)$ are important,  
where $F(D)$ is included in 
the likelihood, $P(\mbox{data}|\mbox{domain}).$ 
To consider our inverse problem in the Bayesian formulation,  
we need the measurability of the function $F.$
Moreover, 
to obtain the stabilities of the posterior and domain ratio  
considered in this study, 
we need the boundedness 
of the image of $\mathcal{D}$ under $F.$
For our model problem, we show such properties using some results 
of \cite{T} that 
studied probabilistic representations of solutions to 
parabolic equations, that is,  
Feynman-Kac type formulas. 
Many reconstruction methods have been considered for 
our model problem and similar problems, e.g., 
\cite{BC2}, \cite{CKY}, 
\cite{HT}, \cite{HNW}, \cite{IK}, 
\cite{K}, and \cite{NW}.
The (locally) 
Lipschitz stability of the posterior density 
ensures that 
the estimation result obtained by such a method with sufficient data
is generally reliable, if
the observational noise is Gaussian and 
$F$ has the above 
measurability and boundedness. 

Although it seems that there are no research results 
for our model problem in the Bayesian formulation, 
many researchers have studied inverse problems 
arising from partial differential equations 
in the Bayesian formulation, e.g., 
\cite{BGh}, \cite{BN},  
\cite{CDS}, \cite{ILS}, \cite{L}, \cite{RSST},    
\cite{vo}, \cite{wangz}, \cite{WMZ}, 
and \cite{Z}. 
In particular, 
for Bayesian geometric inverse problems, 
\cite{ILS} has proposed an 
excellent Bayesian level set method 
with the preconditioned Crank-Nicolson MCMC algorithm.
To apply this method to our model problem, 
we have to establish almost sure continuity of 
the level set map. 
It is an open problem. 
Whether applying the Bayesian level set method 
or the shape identification approach of this paper 
to our model problem with an MCMC algorithm,  
we have to solve (\ref{eq.sect.intro.1})
using some numerical algorithm to calculate the acceptance probability. 
Then, 
if the proposal is accepted, the shape of $\partial C$ is changed
and we have to solve (\ref{eq.sect.intro.1}) on the new domain. 
It has a high computational complexity (cf. \cite{K}). 
Furthermore, if we use our shape identification approach, 
how to select a proposal is also an open problem.  

The remainder of this paper is organized as follows. 
In section \ref{sect.problems},
we state our general shape identification inverse problem,
model problem, and their Bayesian formulations.   
In section \ref{sect.stabilities}, we show the (locally) 
Lipschitz stabilities of the posterior density and domain ratio 
with respect to the observation data 
of our problems. 
Section \ref{sect.appendix} is an appendix. 
In this section, some definitions and a remark 
are described. 

The contents of this paper 
are outlined in Appendix B of \cite{T}. 

\section{Shape identification inverse problems in a Bayesian framework
and model problem} 
\label{sect.problems} 

\subsection{Shape identification inverse problems in a Bayesian framework}
\label{subsect.problems.framework} 

Let $d$ be a positive integer
and denote by $\mathcal{D}^\sharp$ 
the set of bounded domains  (connected open subsets)
in $\mathbb{R}^d.$ 
We consider a finite or infinite subfamily 
$\mathcal{D}$ of $\mathcal{D}^\sharp.$ 
Let $m$ be a positive integer and $F$ be a map from 
$\mathcal{D}$ to $\mathbb{R}^m.$
We call $F$ a \textit{forward operator}. 
We assume that, for a given domain $D \in \mathcal{D},$
we can observe data $y \in \mathbb{R}^m$ 
with noise $\eta \in \mathbb{R}^m,$  
\begin{equation}
y = F(D) + \eta,  
\label{eq.sect.problems.1} 
\end{equation}
where $\eta$ is an $m$-dimensional Gaussian random variable. 
Many researches have considered the case that $\eta$ 
follows an $m$-dimensional Gaussian distribution with some covariance matrix 
or considered 
more general situations (e.g. \cite{DS}). 
For simplicity, we assume that 
the elements $\{\eta_i:\ i = 1, 2, \ldots, m\}$ 
of $\eta = (\eta_i)$ are mutually independent and 
$\eta_i \sim \mathcal{N}(0, \sigma^2),$ 
that is, each $\eta_i$ follows the Gaussian distribution with mean $0$ 
and variance $\sigma^2.$

We are concerned with the following inverse problem. 
\begin{problem}
\label{problem.sect.problems.1} 
Infer $D \in \mathcal{D}$ from given data $y$ of (\ref{eq.sect.problems.1}).  
\end{problem}  
We do not necessarily assume the injectivity of $F.$ 
Therefore, Problem \ref{problem.sect.problems.1} may be 
underdetermined in general. 
In this study, 
we consider Problem \ref{problem.sect.problems.1}  
in a Bayesian framework based on 
\cite{S1} and \cite{DS}. 
In fact, $\mu_0,$ $\Psi,$ $\mu$ 
and $\mu^y$ defined below play the following roles: 
\begin{itemize}
\item
$\mu_0$ is the prior distribution of $D$; 
\item
$\Psi$ is the likelihood; 
\item
$\mu$ is the joint distribution of $(D,y);$ 
\item
$\mu^y$ is the posterior distribution of $D$ given $y$, 
\end{itemize} 
where $\mu^y$ is obtained 
as a disintegration of $\mu.$ 

Let $\mathcal{B}$ be a $\sigma$-field on $\mathcal{D}.$
We always require that $F$ is a measurable map 
on $(\mathcal{D}, \mathcal{B}).$ 
Let $\mu_0$ be a probability measure on $(\mathcal{D}, \mathcal{B})$
and denote by $\|\ \cdot\ \|$ the standard Euclidean norm.
Define 
\begin{equation}
\Psi(D;y) := 
\Frac{1}{(2\pi\sigma^2)^{m/2}}
\exp\left(-\Frac{1}{2\sigma^2}\|y - F(D)\|^2\right)  
\hspace{1cm}
\left(D \in \mathcal{D}, 
y \in \mathbb{R}^m\right)
\label{eq.sect.problems.2} 
\end{equation}
and
\[
Z_\Psi(y) 
:= \int_{\mathcal{D}} \Psi(D;y) \mu_0(dD)
\hspace{1cm}
\left(y \in \mathbb{R}^m\right).
\]
Denote by $dy$ the Lebesgue measure on $\mathbb{R}^m$ 
and define a measure $Q_0$ on $\mathbb{R}^m$ by 
\[
Q_0(dy) := Z_\Psi(y) dy. 
\]
As  
\[
\int_{\mathbb{R}^m} Q_0(dy) 
= \int_{\mathbb{R}^m} \left\{
\int_{\mathcal{D}} \Psi(D;y) \mu_0(dD)\right\} dy
= \int_{\mathcal{D}} \left\{\int_{\mathbb{R}^m} \Psi(D;y) dy\right\} \mu_0(dD)
= 1
\]
from Fubini's theorem,  
$Q_0$ is a probability measure. 
Define a measure $\mu$ on $\mathcal{D} \times \mathbb{R}^m$ 
by
\[
\mu(dD, dy) := \Frac{1}{Z_\Psi(y)} \Psi(D;y) \mu_0(dD) Q_0(dy)
= \Psi(D;y) \mu_0(dD) dy.
\]
Then, $\mu$ is a probability measure, the joint distribution of $(D,y),$ 
and $D$ and $y$ are not independent with respect to $\mu.$
For $y \in \mathbb{R}^m,$ 
define a measure $\mu^y$ on $\mathcal{D}$ by   
\begin{equation}
\mu^y(dD) 
:= \Frac{1}{Z_\Psi(y)} \Psi(D;y) \mu_0(dD). 
\label{eq.sect.problems.3} 
\end{equation}
Then, 
$\mu^y$ is absolutely continuous 
with respect to 
$\mu_0,$ and 
$\mu^y$ satisfies the following from Fubini's theorem:
\begin{enumerate}
\renewcommand{\labelenumi}{\theenumi}
\renewcommand{\theenumi}{(\thesection.\alph{enumi})}
\item
\label{item.sect.problems.1} 
$\mu^y$ is a probability measure on $\mathcal{D}$ for each 
$y \in \mathbb{R}^m;$ 
\item
\label{item.sect.problems.2} 
for every nonnegative measurable function $f$ on   
$\mathcal{D} \times \mathbb{R}^m,$ the function 
\[
y \longmapsto \int_{\mathcal{D}} f(D,y) \mu^y(dD)
\]
is a measurable function on $\mathbb{R}^m;$ 
\item
\label{item.sect.problems.3} 
for every nonnegative measurable function $f$ on   
$\mathcal{D} \times \mathbb{R}^m,$ the equation
\[
\int_{\mathcal{D} \times \mathbb{R}^m} f(D,y) \mu(dD, dy) 
= \int_{\mathbb{R}^m}\left\{\int_{\mathcal{D}} f(D,y) \mu^y(dD)\right\}Q_0(dy)
\]
holds.
\end{enumerate}
Therefore, $\mu^y$ is a disintegration of $\mu$  
(see \cite{CP} p.292 for the definition of disintegration)
and $\mu^y$ is a probability measure 
of the conditional random variable $D|y$ 
from the viewpoint of Kolmogorov's approach to conditioning
(see \cite{Kl} V, Sect. 1, cf. \cite{CP} p.293). 

\begin{lemma}
\label{lemma.sect.problems.1} 
Let $\lambda^y$ be a measure on $\mathcal{D}$ 
indexed by $y \in \mathbb{R}^m.$ 
If $\lambda^y$ is a disintegration of $\mu,$
that is, if $\lambda^y$ 
satisfies \ref{item.sect.problems.1},  
\ref{item.sect.problems.2},  
and \ref{item.sect.problems.3}, 
for every nonnegative measurable function $f$ on   
$\mathcal{D} \times \mathbb{R}^m,$ 
the measurable set  
\[
N_f := \left\{y \in \mathbb{R}^m: 
\ \int_{\mathcal{D}} f(D,y) \mu^y(dD)
\neq \int_{\mathcal{D}} f(D,y) \lambda^y(dD)
\right\}
\]
satisfies $Q_0(N_f) = 0.$   
\end{lemma}  
\textit{Proof.} 
Put 
\begin{eqnarray*}
N_f^+ & := & \left\{y \in \mathbb{R}^m: 
\ \int_{\mathcal{D}} f(D,y) \mu^y(dD)
> \int_{\mathcal{D}} f(D,y) \lambda^y(dD)
\right\},
\\
N_f^- & := & \left\{y \in \mathbb{R}^m: 
\ \int_{\mathcal{D}} f(D,y) \mu^y(dD)
< \int_{\mathcal{D}} f(D,y) \lambda^y(dD)
\right\}.
\end{eqnarray*}
The indicator functions $\chi_{N_f^+}(y)$ and $\chi_{N_f^-}(y)$ 
are measurable. 
If $Q_0(N_f^+) > 0,$ 
we have 
\[
\int_{\mathbb{R}^m}\left\{\int_{\mathcal{D}} 
\chi_{N_f^+}(y) f(D,y) \mu^y(dD)\right\}Q_0(dy)
> 
\int_{\mathbb{R}^m}\left\{\int_{\mathcal{D}} 
\chi_{N_f^+}(y) f(D,y) \lambda^y(dD)\right\}Q_0(dy).
\]
As $\mu^y$ and $\lambda^y$ are disintegrations of $\mu,$
we have 
\begin{eqnarray*}
& & 
\int_{\mathbb{R}^m}\left\{\int_{\mathcal{D}} 
\chi_{N_f^+}(y) f(D,y) \mu^y(dD)\right\}Q_0(dy)
=
\int_{\mathcal{D} \times \mathbb{R}^m} \chi_{N_f^+}(y) f(D,y) \mu(dD, dy) 
\\
& & 
=
\int_{\mathbb{R}^m}\left\{\int_{\mathcal{D}} 
\chi_{N_f^+}(y) f(D,y) \lambda^y(dD)\right\}Q_0(dy).
\end{eqnarray*}
Thus, we have $Q_0(N_f^+) = 0.$
Similarly, we also have $Q_0(N_f^-) = 0.$ 
\hspace*{\fill} 
\QED 
\vspace{0.5cm}

From the general theory as described in \cite{CP}, 
we cannot directly obtain 
the strong uniqueness of the disintegration of $\mu,$  
that is, the property that 
$\cup_f N_f$ is measurable and $Q_0(\cup_f N_f)$ is zero.     
However, $\mu^y$ is a probability measure 
of the conditional random variable $D|y$ 
and we can consider 
(\ref{eq.sect.problems.3}) 
as a Bayesian formula as stated above.
Note that 
$\Psi(D;y)$ and $Z_\Psi(y)$ have the same factor $1/(2\pi\sigma^2)^{m/2}.$ 
Here and in the following, we use the same representation of $\mu^y$ 
as used in \cite{S1} and \cite{DS}, that is, we define
\begin{equation}
\Phi(D;y) := 
\exp\left(-\Frac{1}{2\sigma^2}\|y - F(D)\|^2\right)  
\label{eq.subsect.problems.framework.1} 
\end{equation} 
and 
\begin{equation}
Z(y) 
:= \int_{\mathcal{D}} \Phi(D;y) \mu_0(dD).
\label{eq.sect.problems.Z(y)} 
\end{equation} 
Then, we have
\begin{equation}
\mu^y(dD) 
= \Frac{1}{Z(y)} \Phi(D;y) \mu_0(dD). 
\label{eq.sect.problems.5} 
\end{equation} 
We call $\Phi$ the \textit{potential}.

We also consider, for each $x \in \mathbb{R}^d,$
the conditional probability that the point $x$ is 
included in the unknown domain given the observation data $y.$
Define 
\[
\rho(x|y) := \int_{\mathcal{D}} \chi_D(x) \mu^y(dD)
\hspace{1cm}
\left(x \in \mathbb{R}^d\right),  
\]
where
\[
\chi_D(x) := \left\{
\begin{array}{lll}
1 & & (x \in D)
\\
0 & & (x \in D^c),  
\end{array}
\right.
\]
where $D^c$ is the complement set of $D.$ 
We call $\rho(x|y)$ the \textit{domain ratio} of $x$ under $y.$   

In this study, the conditional probability $\mu^y$ 
and the domain ratio $\rho(\ \cdot\ |y)$ play major roles, 
that is, Problem \ref{problem.sect.problems.1} 
is reformulated as follows. 
\begin{problem}
\label{problem.sect.problems.2} 
Obtain $\mu^y$ on $\mathcal{D}$ 
and $\rho(\ \cdot\ |y)$ on $\mathbb{R}^d$ 
from given data $y$ of (\ref{eq.sect.problems.1}).
\end{problem}

\subsection{Hausdorff distance} 
\label{subsect.problems.hausdorff}

A Lipschitz domain is a domain such that its boundary 
is Lipschitz continuous 
(more precisely, 
see Definition \ref{definition.subsect.problems.framework.1}). 
Our model problem in this study is defined on 
bounded Lipschitz domains 
(see subsection \ref{subsect.problems.model}).
In the following, we denote by $\mathcal{D}_{Lip}$ 
the family of bounded Lipschitz domains 
in $\mathbb{R}^d.$ Then, $\mathcal{D}_{Lip} \subset \mathcal{D}^\sharp.$ 
When we consider our model problem, 
the considered set $\mathcal{D}$ is a  
subfamily of $\mathcal{D}_{Lip},$  
and we assume 
that $\mathcal{D}$ is equipped with 
the Hausdorff distance.
For $\epsilon > 0$ and $A \subset \mathbb{R}^d,$  
put $O_\epsilon(A) := \left\{x \in \mathbb{R}^d: 
\ d(x, A) < \epsilon\right\},$ 
where  
$d(x,A)$ is the Euclidean distance between $x$ and $A.$ 
Then, 
the Hausdorff distance 
between domains $D_1$ and $D_2$ of $\mathcal{D},$
$d_H(D_1, D_2),$ 
is defined by 
\[
d_H(D_1, D_2) := \inf
\left\{\epsilon > 0:\ 
D_1 \subset O_\epsilon(D_2), 
\ D_2 \subset O_\epsilon(D_1) 
\right\}.
\]  
We assume that  
$\mathcal{B}$ is the topological $\sigma$-field 
($\sigma$-field generated by the open sets 
of $\mathcal{D})$ with respect to $d_H,$ and 
require that $F$ is measurable on $(\mathcal{D}, \mathcal{B}).$ 
Note that 
$\mathcal{D}$ is in general \textit{not} a complete metric space if 
$\mathcal{D}$ is an infinite set. 
To apply the results of \cite{S1} and \cite{DS} to our problem directly,  
$\mathcal{D}$ must be complete. 
This completeness ensures that we can use the theory of 
disintegration (cf. \cite{CP}), and then 
the strong uniqueness 
of the disintegration holds.
When $\mathcal{D}$ is not a complete metric space, 
we cannot directly obtain 
the strong uniqueness 
from the general theory as described in \cite{CP}. 
However, we can consider our problem based on \cite{S1} and 
\cite{DS} as mentioned in subsection \ref{subsect.problems.framework}.  
\begin{remark}
\label{remark.subsect.problems.framework.1}
We do not consider the completion of $\mathcal{D}.$ 
The reason is that $F$ can not be extended naturally 
to a function on the completion   
in general.
\end{remark}

\subsection{Model problem} 
\label{subsect.problems.model}

In this study, we consider the following example of 
Problem \ref{problem.sect.problems.2},  
and we call it a \textit{model problem}.
For a subset $U$ of $\mathbb{R}^d,$ 
denote by $\overline{U}$ the closure of $U.$ 
Let $\Omega$ be a bounded Lipschitz domain 
in $\mathbb{R}^d$  
and let $C$ be the union of several (one or more) disjoint simply-connected  
Lipschitz domains (that are not tangential to each other) 
such that $\overline{C} \subset \Omega.$ 
Set $D := \Omega \setminus \overline{C}.$ 
Let $D$ denote a \textit{heat conductor}, and let $C$ denote the union of 
one or more \textit{cavities} inside $D.$   
In the following, we fix $\Omega$ and attempt to identify 
the unknown cavity/cavities 
from all the possible cavities. 
Define $\mathcal{D} = \mathcal{D}_{hc} \subset \mathcal{D}_{Lip}$ 
as the family of such heat conductors with cavities. 
We consider the initial boundary value problem for the heat equation,
(\ref{eq.sect.intro.1}), on $D \in \mathcal{D}_{hc}.$ 
In (\ref{eq.sect.intro.1}),
$T$ is a positive real number, 
$\Delta$ is the standard Laplacian (the coefficient 
$1/2$ is suitable for a Feynman-Kac type formula),
$\nu$ is the unit normal to $\partial \Omega$ directed into the 
exterior of $\Omega,$ and 
$\psi$ is a given $L^2$ function on $[0,T] \times \partial \Omega.$ 

The existence and uniqueness of weak solutions 
to the initial boundary value problem 
(\ref{eq.sect.intro.1})
are well established (see Definition  
\ref{definition.subsect.supple.model.weak} 
for the definition of a weak solution). 
We denote the weak solution by $u^D.$ 
Let $A$ be an open subset of $\partial \Omega$ 
(an accessible portion of $\partial \Omega).$ 
Then, we define a map $G_1$ by 
\[
G_1: \mathcal{D}_{hc} \longrightarrow L^2([0,T] \times A),
\hspace{0.5cm}
D \longmapsto \left.u^D\right|_A. 
\]
Let $G_2$ be a suitable finitization (discretization) map, that is, 
a bounded linear operator  
\[
G_2: L^2([0,T] \times A) \longrightarrow \mathbb{R}^m
\hspace{1cm}
(m \in \mathbb{N}). 
\]
Then, we define a forward operator $F$ by $F := G_2 \circ G_1.$ 
For example, $G_2$ is constructed as follows.
(Other construction methods of $G_2$ have been considered. 
For example, see \cite{ILS}.)   
Discretize $[0,T]$ into $m_T$ intervals,
\[
[0,T] = \mathop{\bigcup}_{i=0}^{m_T-1} [t_i, t_{i+1}],
\]   
where $t_0 = 0,$ $t_{m_T} = T,$ and 
$t_{i+1} - t_i = T/m_T$ $(i = 0, 1, \ldots, m_T-1).$  
Discretize $A$ into $m_A$ subsets, 
\[
A = \mathop{\bigcup}_{j=1}^{m_A} A_j,
\] 
where the Lebesgue measures of $A_1, \ldots, A_{m_A}$ are equal
and the Lebesgue measure of $A_i \cap A_j$ is zero if $i \neq j.$    
Let $m = m_T m_A,$ and define $\{v_{i,j}:\ i=0, 1, \ldots, m_T-1;
j = 1, \ldots, m_A\}$ by 
\begin{equation}
v_{ij} := \int_{t_i}^{t_{i+1}} dt \int_{A_j} v(t,x) dx  
\label{eq.subsect.problems.model.vij}
\end{equation}
for $v \in L^2([0,T] \times A).$ 
This gives 
\[
G_2: L^2([0,T] \times A) \longrightarrow \mathbb{R}^m, 
\hspace{0.5cm}
v \longmapsto \{v_{ij}\}.
\]
Denote by $|A|$ the Lebesgue measure of $A.$ Then, we have 
\begin{equation}
\|G_2(v)\|^2 = \mathop{\sum}_{i=0}^{m_T-1} 
\mathop{\sum}_{j=1}^{m_A} G_2(v)_{ij}^2 
\leq m |A| \|v\|_{L^2([0,T]\times A)}^2.
\label{eq.subsect.problems.model.10}
\end{equation} 
We define the forward operator by 
\begin{equation}
F := G_2 \circ G_1:\  
D \mapsto \left.u^D\right|_A \mapsto 
\left\{\left(\left.u^D\right|_A\right)_{ij}\right\}
\label{eq.subsect.problems.model.F}
\end{equation}
and consider Problem 
\ref{problem.sect.problems.2}. 
To consider this problem, 
we assume that $\mathcal{D} = \mathcal{D}_{hc}$ 
is equipped with 
the Hausdorff distance and 
$\mathcal{B}$ is the 
topological $\sigma$-field.   
As stated in subsection \ref{subsect.problems.hausdorff},
$F$ must be a measurable map on 
$(\mathcal{D},\mathcal{B}) =   
(\mathcal{D}_{hc},\mathcal{B}).$  

For $\alpha \in (0,1],$  
we say that a bounded domain is of class $C^{2,\alpha}$
if its boundary is of class $C^{2,\alpha}$ (more precisely, 
the function $\varphi$ of Definition 
\ref{definition.subsect.problems.framework.1} 
is of class $C^{2,\alpha},$ that is, 
$\varphi$ is of class $C^2$ and all the second-order derivatives 
are $\alpha$-H\"older continuous). 
We consider the following assumption.

\begin{assumption}
\label{assumption.subsect.problems.model.1}
There exists $\alpha \in (0,1]$ such that 
$\Omega$ is of class $C^{2,\alpha}$
and $\psi$ is Lipschitz continuous.  
\end{assumption}
In this study, 
we use a stochastic representation 
(Feynman-Kac type formula) of $u^D,$ 
because this formula gives a direct relation among $u^D,$ 
$\psi,$ and cavities $C$ 
(see (\ref{eq.subsect.formula.2}) and 
(\ref{eq.subsect.problems.model.tau})). 
Under Assumption \ref{assumption.subsect.problems.model.1},
\cite{T} gave the following formula of $u^D,$ 
\begin{equation}
u^D(t,x) = E_x\left[
\int_0^{t\wedge\tau(D)} \psi(t-r, X(r)) L(dr)
\right]
\hspace{1cm}
\left((t,x) \in [0,T] \times \overline{D}\right),
\label{eq.subsect.formula.2}
\end{equation}
where $(X(t), L(t))$ is a pair of stochastic processes. 
The symbols on the right-hand side of 
(\ref{eq.subsect.formula.2})
are defined as follows.

The stochastic process $X$ with $L$ 
is an $\overline{\Omega}$-valued 
diffusion process with a normal reflection on $\partial \Omega$  
starting at $x \in A,$ i.e.,
$X = \{X(t)\}_{0\leq t\leq T}$ is an $\overline{\Omega}$-valued 
continuous stochastic process 
with $X(0) = x,$ and $L = \{L(t)\}_{0\leq t\leq T}$ 
is a continuous, increasing stochastic process 
(local time) such that the pair $(X, L)$ satisfies 
\begin{equation}
\left\{
\begin{array}{ll}
dX(t) = dB(t) + N(X(t)) dL(t)
& \mbox{for $0 \leq t \leq T,$}
\\
L(t) = \displaystyle{\int_0^t 1_{\partial \Omega}}(X(r)) L(dr) 
& \mbox{for $0 \leq t \leq T,$}
\end{array}
\right.
\label{eq.subsect.formula.1}
\end{equation}
where $B$ is a standard Brownian motion on $\mathbb{R}^d$
and $N(X(t))$ is the unit inward normal vector at $X(t) \in \partial \Omega.$   
The existence and uniqueness of such $(X, L),$ 
a strong solution of (\ref{eq.subsect.formula.1}),
have been proved by \cite{Sa} and \cite{Ta}.
Here, $\tau(D)$ is defined by 
\begin{equation}
\tau(D) := \inf \left\{
t:\ 0 \leq t \leq T,\ X(t) \in \overline{C} = \Omega \setminus D 
\right\}
\hspace{1cm}
(\inf \emptyset := \infty),
\label{eq.subsect.problems.model.tau}
\end{equation}
$E_x$ means the expectation considering $X(0) = x,$
and $a \wedge b := \min\{a, b\}.$ 

Then, Theorem 3.1 of \cite{T} implies the following theorem.
\begin{theorem}
\label{theorem.subsect.problems.model.1}
If Assumption \ref{assumption.subsect.problems.model.1} holds,
then for $D \in \mathcal{D}_{hc},$ 
the weak solution $u^D$ of 
(\ref{eq.sect.intro.1})
is represented by (\ref{eq.subsect.formula.2}) and 
$u^D$ is continuous on $[0,T] \times \overline{D}.$ 
\end{theorem}
Theorem 4.1 of \cite{T} 
ensures the continuity of $G_1$ 
with respect to the Hausdorff distance
if Assumption \ref{assumption.subsect.problems.model.1}
is satisfied.
From this fact and (\ref{eq.subsect.problems.model.10}),  
we obtain the following theorem. 
\begin{theorem}
\label{theorem.subsect.problems.model.2}
Under Assumption \ref{assumption.subsect.problems.model.1}, 
the forward operator $F = G_2 \circ G_1$ is continuous 
with respect to the Hausdorff distance on $\mathcal{D}_{hc}.$ 
\end{theorem}
From this theorem, $F$ is a measurable map on 
$(\mathcal{D}, \mathcal{B})
= (\mathcal{D}_{hc}, \mathcal{B})$ 
under Assumption \ref{assumption.subsect.problems.model.1}. 

\section{Stabilities of the posterior density and domain ratio 
with respect to observational data}
\label{sect.stabilities}

\subsection{Stabilities for 
the general shape identification inverse problem}
\label{subsect.stabilities}

Let $\mathcal{D}$ be a subfamily of $\mathcal{D}^\sharp.$
For given data $y$ and $y' \in \mathbb{R}^m,$ 
define  
\[
\sigma(y, y') := \sup\left\{
\|y - F(D)\| \vee \|y' - F(D)\|:\ D \in \mathcal{D}\right\},
\]
where $a \vee b := \max\{a, b\}.$ 
\begin{assumption}
\label{assumption.sect.stabilities.1}
There exists $C_F > 0$ such that $\|F(D)\| < C_F$ 
for every $D \in \mathcal{D}.$ 
\end{assumption}
Note that, under this assumption,
$\sigma(y, y') < \infty$ and   
\[
\min\left(Z(y), Z(y')\right)  
\geq \int_{\mathcal{D}} 
\exp\left(-\Frac{1}{2\sigma^2}\sigma(y,y')\right) 
\mu_0(dD) 
= \exp\left(-\Frac{1}{2\sigma^2}\sigma(y,y')\right) 
> 0. 
\]
The Hellinger distance between $\mu^y$ and $\mu^{y'}$ is defined by 
\begin{eqnarray*}
d_{\mbox{\footnotesize{Hell}}}\left(\mu^y, \mu^{y'}\right)
& := & \sqrt{\Frac{1}{2}
\int_{\mathcal{D}} 
\left\{
\sqrt{\Frac{d\mu^y}{d\mu_0}(D)} 
- 
\sqrt{\Frac{d\mu^{y'}}{d\mu_0}(D)} 
\right\}^2 
\mu_0(dD)}
\\
& = & 
\sqrt{\Frac{1}{2}
\int_{\mathcal{D}} 
\left\{
\Frac{1}{\sqrt{Z(y)}} \sqrt{\Phi(D;y)}
-
\Frac{1}{\sqrt{Z(y')}} \sqrt{\Phi(D;y')}
\right\}^2 
\mu_0(dD)}.
\end{eqnarray*}
We have the following  (locally) Lipschitz 
stabilities of the posterior density and domain ratio 
with respect to the data for our 
general shape identification inverse problem.
Theorem \ref{theorem.sect.stabilities.1}
is similar to Corollary 4.4 of \cite{S1}, 
and Theorem \ref{theorem.sect.stabilities.2}
is a type of stability as stated in 
Remark 4.6 of \cite{DS}. 
\begin{theorem}
\label{theorem.sect.stabilities.1}
If the forward operator $F$ satisfies Assumption 
\ref{assumption.sect.stabilities.1},   
then for given data $y$ and $y' \in \mathbb{R}^m,$ 
\begin{equation}
d_{\mbox{\footnotesize{\rm Hell}}}\left(\mu^y, \mu^{y'}\right)
\leq 
\exp\left(\Frac{3\sigma(y,y')^2}{4\sigma^2}\right)
\Frac{\sigma(y,y')}{\sigma}
\Frac{\|y - y'\|}{\sigma}. 
\label{eq.theorem.sect.stabilities.1.1}
\end{equation}
Furthermore, for every $r > 0,$ there exists $C(r) > 0$ such that  
\[
d_{\mbox{\footnotesize{\rm Hell}}}\left(\mu^y, \mu^{y'}\right)
\leq C(r)\|y - y'\|
\]
for every $y, y'$ with $\|y\| < r, \|y'\| < r.$ 
\end{theorem}
\textit{Proof.} 
We can prove this theorem in the same manner 
as the proof of Theorem 4.2 of \cite{S1}. 
Here, we only prove (\ref{eq.theorem.sect.stabilities.1.1}).
For every $D \in \mathcal{D}$ 
\begin{eqnarray*}
& & 
\left\{
\Frac{1}{\sqrt{Z(y)}} \exp\left(-\Frac{1}{4\sigma^2}\|y - F(D)\|^2\right)  
-
\Frac{1}{\sqrt{Z(y')}} \exp\left(-\Frac{1}{4\sigma^2}\|y' - F(D)\|^2\right)  
\right\}^2 
\\
& \leq & 
\Frac{2}{Z(y)}
\left\{\exp\left(-\Frac{1}{4\sigma^2}\|y - F(D)\|^2\right)  
-
\exp\left(-\Frac{1}{4\sigma^2}\|y' - F(D)\|^2\right)\right\}^2  
\\
& & 
+ 
2 \exp\left(-\Frac{1}{2\sigma^2}\|y' - F(D)\|^2\right)  
\left(\Frac{1}{\sqrt{Z(y)}} - \Frac{1}{\sqrt{Z(y')}}\right)^2. 
\end{eqnarray*}
We estimate the right-hand side. First, 
\begin{eqnarray}
& & 
\left\{\exp\left(-\Frac{1}{4\sigma^2}\|y - F(D)\|^2\right)  
-
\exp\left(-\Frac{1}{4\sigma^2}\|y' - F(D)\|^2\right)\right\}^2  
\label{eq.proof.theorem.sect.stabilities.1.1} 
\\
& \leq & 
\left\{
\Frac{1}{4\sigma^2}\|y - F(D)\|^2
-
\Frac{1}{4\sigma^2}\|y' - F(D)\|^2
\right\}^2 
\nonumber
\\
& = & 
\Frac{1}{4^2\sigma^4}
\left\{(y - y') \iprod (y + y' - 2F(D))\right\}^2 
\nonumber
\\
& \leq & 
\Frac{1}{4^2\sigma^4}
\left\{\|y - y'\| \cdot \|y + y' - 2F(D)\|\right\}^2  
\nonumber
\\
& \leq & 
\Frac{1}{4^2\sigma^4}
\left\{\|y - y'\| \cdot 
(\|y - F(D)\| + \|y' - F(D)\|)\right\}^2  
\nonumber
\\
& \leq & 
\Frac{1}{4} \Frac{\|y-y'\|^2}{\sigma^2} 
\cdot \Frac{\sigma(y,y')^2}{\sigma^2}. 
\nonumber
\end{eqnarray}
Next, 
\[
\left(\Frac{1}{\sqrt{Z(y)}} - \Frac{1}{\sqrt{Z(y')}}\right)^2
=
\Frac{\left(Z(y') - Z(y)\right)^2}{Z(y)Z(y')
\left(\sqrt{Z(y)} + \sqrt{Z(y')}\right)^2},
\]
and similar to (\ref{eq.proof.theorem.sect.stabilities.1.1}),  
\begin{eqnarray}
& & 
\left|Z(y) - Z(y')\right|
\label{eq.proof.theorem.sect.stabilities.1.3} 
\\
& \leq & 
\left|
\exp\left(-\Frac{1}{2\sigma^2}\|y - F(D)\|^2\right) 
-
\exp\left(-\Frac{1}{2\sigma^2}\|y' - F(D)\|^2\right) 
\right|
\nonumber
\\
& \leq & 
\Frac{1}{2\sigma^2}
\|y - y'\| \cdot (\|y - F(D)\| + \|y' - F(D)\|)
\leq
\Frac{\|y-y'\|}{\sigma} 
\cdot \Frac{\sigma(y,y')}{\sigma}.
\nonumber
\end{eqnarray}
Moreover, 
\begin{equation}
\min\{Z(y),Z(y')\} \geq \exp\left(-\Frac{\sigma(y,y')^2}{2\sigma^2}\right) 
\label{eq.proof.theorem.sect.stabilities.1.3.5} 
\end{equation}
and 
\[
Z(y)Z(y')\left(\sqrt{Z(y)} + \sqrt{Z(y')}\right)^2
\geq  
Z(y)Z(y')\left(Z(y) + Z(y')\right)
\geq 2 \exp\left(-\Frac{3\sigma(y,y')^2}{2\sigma^2}\right).  
\]
Therefore, 
\begin{eqnarray*}
& & 
\left\{
\Frac{1}{\sqrt{Z(y)}} \exp\left(-\Frac{1}{4\sigma^2}\|y - F(D)\|^2\right)  
-
\Frac{1}{\sqrt{Z(y')}} \exp\left(-\Frac{1}{4\sigma^2}\|y' - F(D)\|^2\right)  
\right\}^2 
\\
& \leq & 
2 \exp\left(\Frac{3\sigma(y,y')^2}{2\sigma^2}\right)
\Frac{\|y - y'\|^2}{\sigma^2}
\Frac{\sigma(y,y')^2}{\sigma^2}.
\end{eqnarray*}
Thus, we have (\ref{eq.theorem.sect.stabilities.1.1}).
\hspace*{\fill} 
\QED 
\vspace{0.5cm}

\begin{theorem}
\label{theorem.sect.stabilities.2}
If the forward operator $F$ satisfies Assumption 
\ref{assumption.sect.stabilities.1},  
then for given data $y,$ $y' \in \mathbb{R}^m,$ and  
for every $x \in \mathbb{R}^d$ 
we have
\begin{equation}
|\rho(x|y) - \rho(x|y')| 
\leq 
2 \exp\left(\Frac{\sigma(y,y')^2}{2\sigma^2}\right)
\Frac{\sigma(y,y')}{\sigma} 
\Frac{\|y - y'\|}{\sigma}.  
\label{eq.theorem.sect.stabilities.2.1}
\end{equation}
Furthermore, for every $r > 0,$ there exists $C(r) > 0$ such that  
\[
|\rho(x|y) - \rho(x|y')| 
\leq C(r)\|y - y'\|
\]
for every $y, y'$ with $\|y\| < r, \|y'\| < r.$ 
\end{theorem}
\textit{Proof.}
We only prove (\ref{eq.theorem.sect.stabilities.2.1}).
It holds that 
\begin{eqnarray*}
|\rho(x|y) - \rho(x|y')| 
& = & 
\left|\int_{\mathcal{D}} \chi_D(x) \mu^y(dD) 
- 
\int_{\mathcal{D}} \chi_D(x) \mu^{y'}(dD) 
\right|
\\ 
& = & 
\left|\int_{\mathcal{D}} \chi_D(x) \diffrac{\mu^y}{\mu_0}(D) 
\mu_0(dD) 
- 
\int_{\mathcal{D}} \chi_D(x) \diffrac{\mu^{y'}}{\mu_0}(D) 
\mu_0(dD) 
\right|
\\ 
& \leq & 
\int_{\mathcal{D}} \chi_D(x) \left|\diffrac{\mu^y}{\mu_0}(D) 
- \diffrac{\mu^{y'}}{\mu_0}(D) \right|
\mu_0(dD) 
\\ 
& \leq & 
\int_{\mathcal{D}} \left|\diffrac{\mu^y}{\mu_0}(D) 
- \diffrac{\mu^{y'}}{\mu_0}(D) \right|
\mu_0(dD) 
=: 
\left\|\mu^y - \mu^{y'}\right\|_{L^1(\mathcal{D})}.
\end{eqnarray*}
We can estimate this using 
$\left\|\mu^y - \mu^{y'}\right\|_{L^1(\mathcal{D})}
\leq 
2 \sqrt{2}\ d_{\mbox{\footnotesize{Hell}}}\left(\mu^y, \mu^{y'}\right)$ 
and Theorem \ref{theorem.sect.stabilities.1}. 
Here, we estimate $\left\|\mu^y - \mu^{y'}\right\|_{L^1(\mathcal{D})}$ 
directly:  
\begin{eqnarray*}
& & 
\left\|\mu^y - \mu^{y'}\right\|_{L^1(\mathcal{D})}
=
\int_{\mathcal{D}} \left|\diffrac{\mu^y}{\mu_0}(D) 
- \diffrac{\mu^{y'}}{\mu_0}(D) \right|
\mu_0(dD) 
\\
& & 
=  
\int_{\mathcal{D}} 
\left|\Frac{\Phi(D;y)}{Z(y)}
- \Frac{\Phi(D;y')}{Z(y')}
\right| \mu_0(dD) 
\\
& & 
\leq 
\Frac{1}{Z(y)Z(y')} \int_{\mathcal{D}} 
\left\{
\Phi(D;y')|Z(y) - Z(y')|
+ 
Z(y') |\Phi(D;y) - \Phi(D;y')|
\right\} \mu_0(dD) 
\\
& & 
=
\Frac{1}{Z(y)}
\left\{
\left|Z(y) - Z(y')\right|
+ 
\int_{\mathcal{D}}
\left|\Phi(D;y) - \Phi(D;y')\right| \mu_0(dD) 
\right\}
\\
& & 
\leq 
\Frac{2}{Z(y)}
\int_{\mathcal{D}}
\left|\Phi(D;y) - \Phi(D;y')\right| \mu_0(dD). 
\end{eqnarray*}
Then, in the same manner as 
(\ref{eq.proof.theorem.sect.stabilities.1.3}) 
with (\ref{eq.proof.theorem.sect.stabilities.1.3.5}), 
we have 
\begin{equation}
\left\|\mu^y - \mu^{y'}\right\|_{L^1(\mathcal{D})}
\leq 2 \exp\left(\Frac{\sigma(y,y')^2}{2\sigma^2}\right)
\Frac{\sigma(y,y')}{\sigma} \Frac{\|y-y'\|}{\sigma}.  
\end{equation}
Thus, 
we have (\ref{eq.theorem.sect.stabilities.2.1}).
\hspace*{\fill} 
\QED 
\vspace{0.5cm}

As the data $y$ and $y'$ are normally distributed with mean $F(D)$ 
and the components of $y$ and $y'$ have variance $\sigma^2,$
we have 
$E\left[\|y-F(D)\|^2\right]
= E\left[\|y'-F(D)\|^2\right] = m\sigma^2.$ 
Then, 
we probably have $\|y - y'\| \gtrsim \sigma.$ 
In that case, 
Theorems \ref{theorem.sect.stabilities.1} and 
\ref{theorem.sect.stabilities.2} 
are meaningless. 
Therefore, we consider these theorems as follows.
Let $N$ be a sufficiently large integer and 
let $y^{(1)}, y^{(2)}, \ldots, y^{(N)}$ be data. 
Set 
\[
y := \Frac{1}{N} \mathop{\sum}_{i=1}^N y^{(i)}
\] 
and $y' := F(D)$ (the unknown true value). 
Then, $\|y - y'\| \ll \sigma$ and
$\sigma(y,y') \ll \sigma.$ Thus, 
we can consider   
Theorems \ref{theorem.sect.stabilities.1} and 
\ref{theorem.sect.stabilities.2} 
as giving meaningful results 
about the stabilities of $\mu^y$ and 
$\rho(x|y).$ 

\subsection{Case of the model problem} 
\label{subsect.stabilities.model}

Here, we will apply the results of subsection \ref{subsect.stabilities}
to the model problem of subsection \ref{subsect.problems.model}. 

\begin{theorem}
\label{theorem.subsect.stabilities.model.1}
If Assumption \ref{assumption.subsect.problems.model.1} holds, 
then for the forward operator $F$ defined by 
(\ref{eq.subsect.problems.model.F}), 
there exists $C_F > 0$ such that $\|F(D)\| < C_F$ 
for every $D \in \mathcal{D}_{hc}.$ 
Therefore, $F$ satisfies Assumption 
\ref{assumption.sect.stabilities.1}
on $\mathcal{D}_{hc}.$ 
\end{theorem} 
\textit{Proof.} 
Decompose $\psi$ into $\psi(t,x) = \psi_+(t,x) + \psi_-(t,x),$
where 
\[
\psi_+(t,x) := \left\{
\begin{array}{lll} 
\psi(t,x) & & (\psi(t,x) \geq 0)
\\
0 & & (\psi(t,x) < 0),
\end{array}
\right.
\hspace{0.5cm}
\psi_-(t,x) := \left\{
\begin{array}{lll} 
\psi(t,x) & & (\psi(t,x) < 0)
\\
0 & & (\psi(t,x) \geq 0).
\end{array}
\right.
\]
Then, from (\ref{eq.subsect.formula.2}), 
\begin{eqnarray*}
& & 
E_x\left[\int_0^{t\wedge\tau(D)} \psi_-(t-r, X(r)) L(dr)\right]
\\
& & 
\leq u^D(t,x) 
\leq 
E_x\left[\int_0^{t\wedge\tau(D)} \psi_+(t-r, X(r)) L(dr)\right].
\end{eqnarray*}
Therefore, for every $D \in \mathcal{D}_{hc},$  
\[
E_x\left[\int_0^t \psi_-(t-r, X(r)) L(dr)\right]
\leq u^D(t,x) 
\leq 
E_x\left[\int_0^{t} \psi_+(t-r, X(r)) L(dr)\right].
\]
Thus, using (\ref{eq.subsect.problems.model.10}),
we obtain the theorem. 
\hspace*{\fill} 
\QED 
\vspace{0.5cm}

Owing to this theorem,
we can apply the results of subsection \ref{subsect.stabilities}
to the model problem of subsection \ref{subsect.problems.model}
on $\mathcal{D} = \mathcal{D}_{hc},$ 
that is, 
Theorems \ref{theorem.sect.stabilities.1} 
and \ref{theorem.sect.stabilities.2} hold for the problem.  

\section{Appendix. Definitions and remark}
\label{sect.appendix} 

\begin{definition}
[Lipschitz domain] 
\label{definition.subsect.problems.framework.1} 
Let $D$ be a bounded domain of $\mathbb{R}^d.$ 
We call $D$ a bounded Lipschitz domain if, 
for every $x \in \partial D,$ there exist 
a neighborhood $U(x)$ of $x$ and   
a function $\varphi: \mathbb{R}^{d-1} \rightarrow \mathbb{R}$ 
such that:  
\begin{enumerate}
\renewcommand{\labelenumi}{\theenumi}
\renewcommand{\theenumi}{(\roman{enumi})}
\item 
there exists $C_\varphi > 0$ such that
\[
\left|\varphi\left(\xi\right) - 
\varphi\left(\xi'\right)\right| 
\leq C_\varphi\left\|\xi - \xi'\right\|
\] 
for every $\xi, \xi' \in \mathbb{R}^{d-1};$ 
\item
there exists an open neighborhood $O$ of the origin of $\mathbb{R}^d$ 
such that $D \cap U(x)$ can be transformed into  
\[
O \cap \left\{\left(\xi, \xi_d\right) \in \mathbb{R}^d
= \mathbb{R}^{d-1} \times \mathbb{R}: 
\ \xi_d > \varphi\left(\xi\right)\right\}
\]
by a rigid motion $T_x,$ that is, by a rotation plus a translation,
and $T_x(x)$ is the origin.
\end{enumerate}
\end{definition} 

\begin{definition}
[Weak solution to problem (\ref{eq.sect.intro.1})]
\label{definition.subsect.supple.model.weak} 
Put 
\[
\partial_t := \delfrac{}{t}, 
\ \partial_{x_j} := \delfrac{}{x_j} 
\ (j=1, \ldots, d). 
\]
Define
\begin{eqnarray*}
H^1(D) & := & 
\left\{f \in L^2(D):\ \partial_{x_j} f \in L^2(D)
\ (j=1, \ldots, d)
\right\},
\\
C([0,T]; L^2(D)) & := & 
\left\{f: [0,T] \rightarrow L^2(D):\ \mbox{$f$ is continuous}\right\},
\\
V^{0,1}((0,T) \times D) & := &  
C([0,T]; L^2(D)) \cap L^2((0,T); H^1(D)),
\\
H^{1,1}((0,T) \times D) & := & 
\left\{f \in L^2((0,T) \times D):\ 
\partial_t f, \partial_{x_j} f \in L^2((0,T) \times D)
\ (j=1, \ldots, d) \right\},
\end{eqnarray*}
and
\[
\nabla_x f := \left(\partial_{x_1} f, \ldots, \partial_{x_d} f\right). 
\]
The function spaces $H^1(D)$ and $H^{1,1}((0,T) \times D)$ are equipped with 
the usual Sobolev norms, and 
the function space $V^{0,1}((0,T) \times D)$ is equipped with 
the norm $\| \cdot \|_{D}$ defined by 
\[
\|u\|_D := \mathop{\max}_{0 \leq t \leq T} 
\|u(t,\cdot)\|_{L^2(D)} +  \|\nabla_x u\|_{L^2((0,T) \times D)}. 
\]
Denote by $\gamma$ the boundary trace operator from $H^1(D)$  
into $L^2(\partial D).$ 
A \textit{weak solution} $u^D$ to problem (\ref{eq.sect.intro.1}) 
is a function such that it belongs to $V^{0,1}((0,T) \times D)$ and 
satisfies the weak form of (\ref{eq.sect.intro.1}), 
that is: 
\begin{enumerate}
\renewcommand{\labelenumi}{\theenumi}
\renewcommand{\theenumi}{(\roman{enumi})}
\item
$\left.u^D\right|_{\partial C} = 0,$ 
that is, $\gamma u^D(t,\ \cdot\ ) = 0$ on $\partial C$ 
for almost every $t \in (0,T);$ 
\item
for every $\eta \in H^{1,1}(D)$ with 
$\left.\eta\right|_{\partial C} = 0$ and $\eta(T,\ \cdot\ ) = 0,$ 
it holds that 
\[
\int_D u^D \partial_t \eta\ dtdx
- \int_D \nabla_x u^D \iprod \nabla_x \eta\ dtdx  
+ \int_{\partial \Omega} \psi \gamma\eta\ dt S(dx) = 0,   
\]
where $S(dx)$ is the Lebesgue measure on $\partial \Omega.$
\end{enumerate}
\end{definition} 

\begin{remark}
[Stochastic representation formula for a general parabolic problem] 
\label{remark.sect.appendix.st.repre} 
A rather more general (backward) parabolic boundary value problem than
(\ref{eq.sect.intro.1}) was considered in \cite{T}.
In particular, the Dirichlet boundary and the Neumann boundary 
are allowed to meet, and the Dirichlet boundary is allowed 
to vary with time.   
For the weak solution to the problem, the stochastic representation
and continuity property, which are generalizations of 
Theorems \ref{theorem.subsect.problems.model.1} 
and \ref{theorem.subsect.problems.model.2},  
were shown 
through the coupled
martingale formulation for the basic diffusion process.
If the weak solution $u^D$ to (\ref{eq.sect.intro.1})
has a smooth extension belonging to
$C^{1,2}\left([0,T] \times \mathbb{R}^d\right),$
we can prove Theorem \ref{theorem.subsect.problems.model.1}
using It\^o's formula.
However, in a general case,  
we have to apply the formula to the weak solution via a suitable 
smooth approximation
procedure for the solution as in \cite{T}.
\end{remark} 

\section*{Acknowledgements}
The author is most grateful to 
Professor Masaaki Tsuchiya for his valuable discussion 
and advice.
The author also appreciates the helpful comments from    
the anonymous reviewers. 

This research did not receive any specific grant from funding 
agencies in the public, commercial, or not-for-profit sectors.

\bibliographystyle{elsarticle-num} 

\begin{thebibliography}{00}
\bibitem{BGh}
B.-Thanh, T. and Ghattas, O.,
An analysis of infinite dimensional Bayesian inverse shape 
acoustic scattering and its numerical approximation,
\textsl{SIAM/ASA Journal on Uncertainty Quantification},
2, 203-222 (2014)
DOI: 10.1137/120894877
\bibitem{BN}
B.-Thanh, T. and Nguyen, Q. P.,
FEM-based discretization-invariant MCMC methods for
PDE-constrained Bayesian inverse problems,  
\textsl{Inverse Problems and Imaging},
10(4), 943-975 (2016)
DOI: 10.3934/ipi.2016028
\bibitem{BCSV}
Bacchelli, V., Di Cristo, M., Sincich, E., and Vessella, S.,
A parabolic inverse problem with mixed
boundary data. Stability estimates for the
unknown boundary and impedance, 
\textsl{Transactions of the American Mathematical Society},
366(8), 3965–3995 (2014)
DOI: 10.1090/S0002-9947-2014-05807-8 
\bibitem{BC1}
Bryan, K. and Caudill, L. F.,
Uniqueness for a boundary identification problem in thermal imaging,
\textsl{Electronic Journal of Differential Equations},
Conference 01 C(1), 23-39 (1997)
\bibitem{BC2}
Bryan, K. and Caudill, L., 
Reconstruction of an unknown boundary portion from Cauchy data in $n$ 
dimensions, \textsl{Inverse Problems}, 
21, 239-255 (2005)  
DOI: 10.1088/0266-5611/21/1/015
\bibitem{CP} 
Chang, J. T. and Pollard, D.,
Conditioning as disintegration, 
\textsl{Statistica Neerlandica}, 
51(3), 287-317 (1997) 
DOI: 10.1111/1467-9574.00056
\bibitem{CKY}
Chapko, R. and Kress, R. and Yoon, J. R.,
An inverse boundary value problem for the heat equation: The
Neumann condition,
\textsl{Inverse Problems},  
15, 1033–1049 (1999)
DOI: 10.1088/0266-5611/15/4/313
\bibitem{CDS}
Cotter, S. L., Dashti, M., and Stuart, A. M., 
Approximation of Bayesian Inverse Problems, 
\textsl{SIAM Journal of Numerical Analysis},
48(1), 322-345 (2010)
DOI: 10.1137/090770734
\bibitem{CRV}
Cristo, M. D., Rondi, L., and Vessella, S.,
Stability properties of an inverse parabolic problem 
with unknown boundaries,
\textsl{Annali di Matematica Pura ed Applicata}, 185(2), 223-255 (2006) 
DOI: 10.1007/s10231-005-0152-x
\bibitem{DS}
Dashti, M. and Stuart, A. M.,
The Bayesian approach to inverse problems,
\textsl{Handbook of Uncertainty Quantification}, 
1-118 (2016)
DOI: 10.1007/978-3-319-11259-6\_7-1 
\bibitem{HT}
Harbrecht, H. and Tausch, J.,
An efficient numerical method for a shape-identification 
problem arising from the heat equation,
\textsl{Inverse Problems}, 27, 065013 (2011)  
DOI: 10.1088/0266-5611/27/6/065013
\bibitem{HNW}
Heck, H., Nakamura, G., and Wang, H.,
Linear sampling method for identifying cavities in a heat conductor,
\textsl{Inverse Problems}, 
29, 075014 (2012) 
DOI: 10.1088/0266-5611/28/7/075014
\bibitem{ILS}
Iglesias, M. A., Lu, Y., and Stuart, A. M.,
A Bayesian level set method for geometric inverse problems, 
\textsl{Interfaces and Free Boundaries}, 
18(2), 181-217 (2016)
DOI: 10.4171/IFB/362 
\bibitem{IK}
Ikehata, M. and Kawashita, M.,
The enclosure method for the heat equation, 
\textsl{Inverse Problems}, 
25, 075005 (2009)  
DOI: 10.1088/0266-5611/25/7/075005
\bibitem{K}
Kawakami, H., 
Reconstruction algorithm for unknown cavities 
via Feynman-Kac type formula, 
\textsl{Computational Optimization and Applications}, 
61(1), 101-133 (2015)
DOI: 10.1007/s10589-014-9706-4
\bibitem{KT}
Kawakami, H. and Tsuchiya, M.,
Uniqueness in shape identification of a time-varying domain
and related parabolic equations on non-cylindrical domains,  
\textsl{Inverse Problems}, 
26, 125007 (2010)
DOI: 10.1088/0266-5611/26/12/125007
\bibitem{Kl}
Kolmogorov, A. N., 
\textsl{Foundation of the Theory of Probability
(English translation)}, Chelsea Pub. Co., New York (1956) 
\bibitem{L}
Litvinenko, A., 
Partial inversion of elliptic operator to speed up
computation of likelihood in Bayesian inference,
\textsl{arXiv}, 1708.02207 v1 (2017) 
\bibitem{NW}
Nakamura, G. and Wang, H.,
Reconstruction of an unknown cavity with Robin boundary 
condition inside a heat conductor,  
\textsl{Inverse Problems}, 
31(12), 125001 (2015) 
DOI: 10.1088/0266-5611/31/12/125001
\bibitem{RSST}
Ruggeri, F., Sawlan, Z., Scavino, M., and Tempone, R., 
A Hierarchical Bayesian Setting for an Inverse
Problem in Linear Parabolic PDEs
with Noisy Boundary Conditions, 
\textsl{Bayesian Analysis}, 
12(2), 407-433 (2017) 
DOI: 10.1214/16-BA1007
\bibitem{Sa}
Saisho, S., 
Stochastic differential equations for multi-dimensional domain with 
reflecting boundary,
\textsl{Probability Theory and Related Fields},  
74, 455-477 (1987)
DOI: 10.1007/BF00699100
\bibitem{S1}
Stuart, A. M.,
Inverse Problems: a Bayesian perspective,
\textsl{Acta Numerica}, 
19, 451-559 (2010)   
DOI: 10.1017/S0962492910000061
\bibitem{Ta}
Tanaka, H., 
Stochastic differential equations with reflecting 
boundary condition in convex regions, 
\textsl{Hiroshima Mathematical Journal},
9, 163-177 (1979)
\bibitem{T}
Tsuchiya, M.,
Probabilistic representation of weak solutions to a 
parabolic equation with a mixed boundary condition 
on a non smooth domain 
(Appendix B by Kawakami, H.),
\textsl{arXiv}, 1710.05136 v1 (2017)
\bibitem{vo}
Vollmer, S. J.,
Posterior consistency for Bayesian inverse
problems through stability and regression results,      
\textsl{Inverse Problems}, 
29(12), 125011 (2013)
DOI: 10.1088/0266-5611/29/12/125011
\bibitem{wangz}
Wang, J. and Zabaras, N.,
Using Bayesian statistics in the estimation of heat 
source in radiation,
\textsl{International Journal of Heat and Mass Transfer}, 
48, 15-29 (2005)
DOI: 10.1016/j.ijheatmasstransfer.2004.08.009
\bibitem{WMZ}
Wang, Y., Ma, F., and Zheng, E., 
Bayesian method for shape reconstruction in the inverse
interior scattering problem,
\textsl{Mathematical Problems in Engineering},  
935294 (2015)
DOI: 10.1155/2015/935294
\bibitem{Z}
Zambelli, A. E., 
A Multiple Prior Monte Carlo Method for
the Backward Heat Diffusion Problem,
\textsl{Proceedings of the 11th International Conference
on Computational and Mathematical Methods
in Science and Engineering, CMMSE 2011},   
26-30 (2011) 
\end{thebibliography}


\end{document}


\newtheorem{assumption}{Assumption}[section]
\newtheorem{lemma}{Lemma}[section]
\newtheorem{theorem}{Theorem}[section]
\newtheorem{problem}{Problem}[section]
\newtheorem{remark}{Remark}[section]
\newtheorem{definition}{Definition}[section]
\newcommand{\Frac}[2]{\displaystyle\frac{#1}{#2}}
\newcommand{\delfrac}[2]{\displaystyle\frac{\partial #1}{\partial #2}}
\newcommand{\iprod}{\, \mbox{\raise.3ex\hbox{\tiny $\bullet$}}\, }
\newcommand{\diffrac}[2]{\displaystyle\frac{d #1}{d #2}}
\newcommand{\QED}{$\Box$} 

\begin {center}
\begin{minipage}{14cm}  
\textbf{\Large{Uniqueness of disintegration: Addendum to 
``Stabilities of shape identification inverse problems in 
a Bayesian framework''
[J. Math. Anal. Appl. (2020) 123903]}} 
\end{minipage}
\\
[0.5cm]
Hajime Kawakami
\\
[0.2cm]
Mathematical Science Course,
Akita University, Akita 010-8502, Japan
\\
[0.2cm]
kawakami@math.akita-u.ac.jp
\end{center}
\vspace{0.5cm} 

\noindent
\textbf{Abstract}
\vspace{0.2cm} 

\noindent
The strong uniqueness of disintegration holds 
if the $\sigma$-algebra on the considered space
is countably generated. 
This strong uniqueness yields 
the strong uniqueness of the posterior
as a disintegrating measure in 
the Bayesian inverse problem
discussed in \cite{HKss}.  
\vspace{0.5cm}

\noindent
\textit{keywords}:
\hspace{0.1cm} 
Disintegration,\ Uniqueness,  
\ Bayesian inverse problem,\ Posterior 

\noindent
\textit{2010 MSC}:
\hspace{0.1cm} 
28A50,\ 35R30,\ 60A10,\ 62F15

\section{Introduction} 

In this addendum to \cite{HKss}, we consider
the weak and strong uniqueness 
of disintegration 
and an application of strong uniqueness 
to the Bayesian inverse problem
discussed in \cite{HKss}.
Experts will already 
know or be able to show the results of this addendum, 
but the author believes that the results are worth stating 
for a more general audience. 

As stated in \cite{CP} on p.294, we can consider 
a conditional probability as a disintegrating measure. 
The weak uniqueness of disintegrations yields 
the well-definedness of 
conditional probability measures 
from the viewpoint of Kolmogorov's approach to conditioning
(see \cite{Kl} V, Sect. 1, cf. \cite{CP} p.293).
In \cite{HKss}, the author considered  
shape identification inverse problems in a Bayesian framework.
The posterior is considered 
as a disintegration of the 
joint distribution.
In \cite{HKss}, 
this posterior is constructed specifically 
(therefore, the existence of the posterior is ensured)
and only weak uniqueness is proved.
Thus, strong uniqueness is the remaining problem.

According to Theorem 1 (existence theorem) in \cite{CP}, 
a $\sigma$-finite Radon measure on a metric space has
a disintegration and the strong uniqueness of that  
disintegration holds. 
It is unclear whether or not the
joint distribution in \cite{HKss}  
is a Radon measure. 
However, it can be shown that   
the strong uniqueness of disintegration holds generally if 
the $\sigma$-algebra on a considered space is countably generated
(see Theorem \ref{theorem.sect.uniqueness.1}).
It can be also shown that   
the $\sigma$-algebra on the considered (product) space in \cite{HKss}
is countably generated 
(see Lemma \ref{lemma.sect.uniqueness.1} and the proof of
Theorem \ref{theorem.sect.uniqueness.appl.1}). 
Therefore, the posterior in \cite{HKss} satisfies 
strong uniqueness as a disintegrating measure 
(see Theorem \ref{theorem.sect.uniqueness.appl.1}). 

\section{Uniqueness theorem of disintegration}
\label{sect.uniqueness} 

Let $\left(\mathcal{X}, \mathcal{A}_{\mathcal{X}}, \mu_{\mathcal{X}}\right)$
and $\left(\mathcal{Y}, \mathcal{A}_{\mathcal{Y}}, \mu_{\mathcal{Y}}\right)$
be measure spaces, where 
$\mathcal{A}_{\mathcal{X}}$ and $\mathcal{A}_{\mathcal{Y}}$
are $\sigma$-algebras, 
and $\mu_{\mathcal{X}}$ and $\mu_{\mathcal{Y}}$ are $\sigma$-finite measures.   
Let $T: \left(\mathcal{X}, \mathcal{A}_{\mathcal{X}}\right) \rightarrow 
\left(\mathcal{Y}, \mathcal{A}_{\mathcal{Y}}\right)$
be a measurable map.
Assume that 
$\{y\} \in \mathcal{A}_{\mathcal{Y}}$ for every $y \in \mathcal{Y}.$ 
Let $\mathcal{M}$ denote the set of 
nonnegative measurable functions on $\mathcal{X}.$  
Then, a family of measures on $\mathcal{A}_{\mathcal{X}},$
$\{\mu^y\}$ $(y \in \mathcal{Y}),$ 
is called a $(T,\mu_{\mathcal{Y}})$-disintegration
of $\mu_{\mathcal{X}}$ if the following
properties hold (cf. Definition 1 in \cite{CP}):
\begin{enumerate}
\renewcommand{\labelenumi}{\theenumi}
\renewcommand{\theenumi}{(\roman{enumi})}
\item
\label{item.sect.uniqueness.i} 
$\mu^y$ is a $\sigma$-finite measure on $\mathcal{A}_{\mathcal{X}}$
and $\mu^y\left(\mathcal{X} \setminus T^{-1}(y)\right) = 0$ 
for $\mu_{\mathcal{Y}}$-almost every $y;$
\end{enumerate}
and, for each $f \in \mathcal{M},$ 
\begin{enumerate}
\renewcommand{\labelenumi}{\theenumi}
\renewcommand{\theenumi}{(\roman{enumi})}
\setcounter{enumi}{1} 
\item
\label{item.sect.uniqueness.ii}
$y \mapsto \displaystyle{\int_{\mathcal{X}}} f(x) \mu^y(dx)$
is a measurable function; 
\item
\label{item.sect.uniqueness.iii}
$\displaystyle{\int_{\mathcal{X}}} f(x) \mu_{\mathcal{X}}(dx)
= \displaystyle{\int_{\mathcal{Y}}}\left\{
\displaystyle{\int_{\mathcal{X}}} f(x) \mu^y(dx)
\right\}\mu_{\mathcal{Y}}(dy).$
\end{enumerate}
Let $\mu_1^y$ and $\mu_2^y$ be
$(T,\mu_{\mathcal{Y}})$-disintegrations of $\mu_{\mathcal{X}}.$ 
For each $f \in \mathcal{M},$ 
we set 
\[
N_f :=  
\left\{y \in \mathcal{Y}: 
\displaystyle{\int_{\mathcal{X}}} f(x) \mu_1^y(dx)
\neq 
\displaystyle{\int_{\mathcal{X}}} f(x) \mu_2^y(dx)
\right\}.
\]
According to \ref{item.sect.uniqueness.ii},
$N_f$ is a measurable set. 
We define the weak uniqueness and strong uniqueness of
the $(T,\mu_{\mathcal{Y}})$-disintegration of $\mu_{\mathcal{X}}$ as follows: 
\begin{itemize}
\item
Weak uniqueness implies that 
\[
\mu_{\mathcal{Y}}\left(N_f\right) = 0
\]
holds for every $f \in \mathcal{M}$  
and every pair $(\mu_1^y, \mu_2^y)$ of  
$(T,\mu_{\mathcal{Y}})$-disintegrations of $\mu_{\mathcal{X}}.$ 
\item
Strong uniqueness implies that $\bigcup_{f \in \mathcal{M}} N_f$
is a measurable set and   
\[
\mu_{\mathcal{Y}}\left(\bigcup_{f \in \mathcal{M}} N_f \right) = 0
\]
holds for every pair $(\mu_1^y, \mu_2^y)$ of  
$(T,\mu_{\mathcal{Y}})$-disintegrations of $\mu_{\mathcal{X}}.$ 
\end{itemize}
\begin{theorem}
\label{theorem.sect.uniqueness.1} 
The weak uniqueness of the   
$(T,\mu_{\mathcal{Y}})$-disintegration of $\mu_{\mathcal{X}}$ always holds. 
If $\mathcal{A}_{\mathcal{X}}$ is countably generated (i.e., 
there exists a countable subset $\mathcal{C}$ of $\mathcal{A}_{\mathcal{X}}$ 
such that $\sigma(\mathcal{C}) = \mathcal{A}_{\mathcal{X}}),$ 
then the strong uniqueness of the  
$(T,\mu_{\mathcal{Y}})$-disintegration of $\mu_{\mathcal{X}}$ also holds. 
\end{theorem}
\textit{Proof.}
First, we prove weak uniqueness.
Let $\mu_1^y$ and $\mu_2^y$ be
$(T,\mu_{\mathcal{Y}})$-disintegrations of $\mu_{\mathcal{X}}.$ 
For each $f \in \mathcal{M},$ 
we set  
\begin{eqnarray*}
N_f^+ 
& := & 
\left\{y \in \mathcal{Y}: 
\displaystyle{\int_{\mathcal{X}}} f(x) \mu_1^y(dx)
> 
\displaystyle{\int_{\mathcal{X}}} f(x) \mu_2^y(dx)
\right\},
\\
N_f^-
& := & \left\{y \in \mathcal{Y}: 
\displaystyle{\int_{\mathcal{X}}} f(x) \mu_1^y(dx)
<  
\displaystyle{\int_{\mathcal{X}}} f(x) \mu_2^y(dx)
\right\}.
\end{eqnarray*}
From \ref{item.sect.uniqueness.ii},  
$N_f^+$ and $N_f^-$ are measurable sets. 
Let $\chi_A$ denote 
the indicator function of a set $A.$
From \ref{item.sect.uniqueness.i},  
\[
N_f^+ = 
\left\{y \in \mathcal{Y}: 
\displaystyle{\int_{\mathcal{X}}} \chi_{T^{-1}(\{y\})}(x)
f(x) \mu_1^y(dx)
> 
\displaystyle{\int_{\mathcal{X}}} \chi_{T^{-1}(\{y\})}(x)
f(x) \mu_2^y(dx)
\right\}.
\]
Note that for $j = 1, 2$ and
almost every $y \in N_f^+,$
\[
\displaystyle{\int_{\mathcal{X}}} \chi_{T^{-1}(\{y\})}(x)
f(x) \mu_j^y(dx)
=
\displaystyle{\int_{\mathcal{X}}} \chi_{T^{-1}(N_f^+)}(x)
f(x) \mu_j^y(dx).
\]
Therefore, if $\mu_{\mathcal{Y}}\left(N_f^+\right) \neq 0,$ 
then 
\begin{eqnarray*}
& & 
\int_{\mathcal{Y}}
\left\{ 
\displaystyle{\int_{\mathcal{X}}} \chi_{T^{-1}(N_f^+)}(x)
f(x) \mu_1^y(dx)
\right\}
\mu_{\mathcal{Y}}(dy)
\\
& > & 
\int_{\mathcal{Y}}
\left\{ 
\displaystyle{\int_{\mathcal{X}}} \chi_{T^{-1}(N_f^+)}(x)
f(x) \mu_2^y(dx)
\right\}
\mu_{\mathcal{Y}}(dy).
\end{eqnarray*}
In contrast, from 
\ref{item.sect.uniqueness.iii},    
\begin{eqnarray*}
& & 
\int_{\mathcal{Y}}
\left\{ 
\displaystyle{\int_{\mathcal{X}}} \chi_{T^{-1}(N_f^+)}(x)
f(x) \mu_1^y(dx)
\right\}
\mu_{\mathcal{Y}}(dy)
\\
& = & 
\int_{\mathcal{Y}}
\left\{ 
\displaystyle{\int_{\mathcal{X}}} \chi_{T^{-1}(N_f^+)}(x)
f(x) \mu_2^y(dx)
\right\}
\mu_{\mathcal{Y}}(dy).
\end{eqnarray*}
Therefore, $\mu_{\mathcal{Y}}\left(N_f^+\right) = 0.$ Similarly,
$\mu_{\mathcal{Y}}\left(N_f^-\right)  = 0.$
Finally, we have $\mu_{\mathcal{Y}}\left(N_f\right)  = 0,$ 
which proves weak uniqueness.

Next, we prove strong uniqueness.
It is equivalent to that 
\begin{equation}
\mu_{\mathcal{Y}}
\left(\left\{y \in \mathcal{Y}:\ \mu_1^y \neq \mu_2^y\right\}\right) = 0
\label{eq.pf.theorem.sect.uniqueness.1.1} 
\end{equation}
holds if $\mu_1^y$ and $\mu_2^y$ are
$(T,\mu_{\mathcal{Y}})$-disintegrations of $\mu_{\mathcal{X}}.$ 
Let $\widetilde{\mathcal{A}}_{\mathcal{X}}$ 
be the algebra 
generated by $\mathcal{C}.$ Then,  
$\mathcal{A}_{\mathcal{X}}
= \sigma\left(\widetilde{\mathcal{A}}_{\mathcal{X}}\right)$ and  
$\widetilde{\mathcal{A}}_{\mathcal{X}}$ is countable 
because the family of all finite subsets of a countable set 
is also countable. 
We set 
\[
N := \left\{ y \in \mathcal{Y}:\ 
\exists A \in \widetilde{\mathcal{A}}_{\mathcal{X}}, 
\mu_1^y(A) \neq \mu_2^y(A)
\right\}.  
\]
Then, based on weak uniqueness, we have 
$\mu_{\mathcal{Y}}(N) = 0.$ 
For every $y \not\in N$ and every  
$A \in \widetilde{\mathcal{A}}_{\mathcal{X}},$ 
\[
\mu_1^y(A) = \mu_2^y(A). 
\]
Because $\mu_1^y$ and $\mu_2^y$ are $\sigma$-finite and 
$\mathcal{A}_{\mathcal{X}}
= \sigma\left(\widetilde{\mathcal{A}}_{\mathcal{X}}\right),$ 
we have 
\[
\mu_1^y(A) = \mu_2^y(A)
\]
for every $y \not\in N$ and  
every $A \in \mathcal{A}_{\mathcal{X}}$  
according to Dynkin's extension theorem (Lemma 1.1 in \cite{Dy}, 
Theorem 8.5 in \cite{A})  
or the extension theorem called  
Hahn-Kolmogorov's theorem 
(Theorem 1.7.8 in \cite{T}), 
Carath\'eodory's theorem (Theorem 6.2 in \cite{A}), 
E. Hopf's theorem (5.5 in \cite{LM}),   
or other names (cf. 5.16 in \cite{LM}).   
Therefore, we have (\ref{eq.pf.theorem.sect.uniqueness.1.1}),  
which proves strong uniqueness.
\hspace*{\fill} 
\QED

\begin{remark}
In \cite{CP}, the strong uniqueness of the disintegration of 
a $\sigma$-finite Radon measure on a metric space $\mathcal{X}$ is 
proved by reducing it to the case where
$\mathcal{X}$ is a compact metric space.   
Because the Borel field 
(i.e., the $\sigma$-algebra generated by the open sets) 
on a compact metric space is countably generated, 
this result also follows from Theorem \ref{theorem.sect.uniqueness.1}.   
\end{remark} 

\section{Application of Theorem \ref{theorem.sect.uniqueness.1} 
to the framework in \cite{HKss}} 

From Subsection 2.2 and onward in \cite{HKss}, 
the disintegration is considered in the following context.
Let $d$ be a positive integer and 
$\mathcal{D}^\sharp$ 
be the family of all bounded domains 
(i.e., bounded connected open subsets)
in $\mathbb{R}^d.$ 
Assume that $\mathcal{D}^\sharp$ is equipped with 
the Hausdorff distance $d_H.$   
Let $\mathcal{D}$ be a metric subspace of $\mathcal{D}^\sharp,$ 
$\mathcal{B}_{\mathcal{D}}$ be the Borel field 
of $\mathcal{D},$ 
and $\mathcal{B}_{\mathbb{R}^m}$ be the Borel field of $\mathbb{R}^m,$
where $m$ is a positive integer. 
Assume that 
$\mathcal{X} = \mathcal{D} \times \mathbb{R}^m$
(the product metric space),  
$\mathcal{A}_{\mathcal{X}}
= \mathcal{B}_{\mathcal{D}} \times \mathcal{B}_{\mathbb{R}^m}$ 
(the product $\sigma$-algebra),  
$\mathcal{Y} = \mathbb{R}^m,$ 
$\mathcal{A}_{\mathcal{Y}} = \mathcal{B}_{\mathbb{R}^m},$
$T: \mathcal{X} \rightarrow \mathcal{Y}$ is defined by $T(D,y) = y$  
$((D,y) \in \mathcal{D} \times \mathbb{R}^m),$ 
and $\mu_{\mathcal{X}}$ and 
$\mu_{\mathcal{Y}}$ are probability measures
($\mu_{\mathcal{X}}$ is a joint distribution).
Then, 
a $(T,\mu_{\mathcal{Y}})$-disintegration
$\mu^y$ of $\mu_{\mathcal{X}}$ 
is the posterior, that is, 
the probability measure
of the conditional random variable $D|y,$
and we have 
\begin{eqnarray*}
& & 
\int_{\mathcal{X}} f(D,y)\mu^y(dD)
= \int_{T^{-1}(y)} f(D,y)\mu^y(dD)
\\
& & 
= \int_{\mathcal{D}\times\{y\}} f(D,y)\mu^y(dD)
= \int_{\mathcal{D}} f(D,y)\mu^y(dD)
\end{eqnarray*}
for $\mu_{\mathcal{Y}}$-almost every $y$ and 
every $f \in \mathcal{M}.$

In \cite{HKss}, the author only proved the weak uniqueness 
of the disintegration in this context.
However, the following theorem holds.
\begin{theorem}
\label{theorem.sect.uniqueness.appl.1}
In the context described above, 
the strong uniqueness of the   
$(T,\mu_{\mathcal{Y}})$-disintegration of $\mu_{\mathcal{X}}$ holds.
\end{theorem}
This theorem is derived from the following 
Lemma \ref{lemma.sect.uniqueness.1} 
and Theorem \ref{theorem.sect.uniqueness.1}. 
\begin{lemma}
\label{lemma.sect.uniqueness.1} 
Every subspace $\mathcal{D}$ of $\mathcal{D}^\sharp$ is a separable
metric space. 
\end{lemma}
\textit{Proof.}
If $\mathcal{D}^\sharp$ is separable, 
then every subspace $\mathcal{D}$ is also separable.
We will now prove the separability of $\mathcal{D}^\sharp.$ 
For $a = \left(a_1, \ldots, a_d\right) \in \mathbb{Z}^d$ 
and $\epsilon > 0,$ we set  
\[
Q(a,\epsilon) 
:= [a_1\epsilon, (a_1+1)\epsilon]
\times \cdots \times [a_d\epsilon, (a_d+1)\epsilon].
\]  
For $D \in \mathcal{D}^\sharp,$ we set 
\[
\mathcal{Q}(D, \epsilon) 
:= \left\{Q(a,\epsilon):\ a \in \mathbb{Z}^d, 
\ Q(a,\epsilon) \cap D \neq \emptyset
\right\}.
\]
Then, $\mathcal{Q}(D, \epsilon)$ is a finite set for every $\epsilon > 0.$ 
We set 
\[
S(D,\epsilon) :=
\mbox{Int}\left(\mathop{\bigcup}_{Q(a,\epsilon) \in \mathcal{Q}(D, \epsilon)}
Q(a,\epsilon)\right),
\]
where ``Int'' indicates the interior.
Then, we have 
$d_H(D, S(D,\epsilon)) < \sqrt{d} \epsilon.$
Note that 
\[
\mathcal{Q} 
:= 
\left\{Q\left(a,\Frac{1}{n}\right):\ a \in \mathbb{Z}^d,
\ n \in \mathbb{N}  
\right\}
\]
is a countable set. Therefore, 
$\mathcal{Q}$ can be rewritten as    
$\mathcal{Q} = \left\{ Q_m:\ m \in \mathbb{N}\right\}.$ 
We define
\begin{eqnarray*}
\mathcal{Q}^\ast := \left\{
\mbox{Int}\left(\mathop{\bigcup}_{m \in I} Q_m\right)
\right.
& : & I \subset \mathbb{N}, 
\ |I| < \infty,  
\ Q_m \in \mathcal{Q}, 
\\
& & 
\left.
\mbox{Int}\left(\mathop{\bigcup}_{m \in I} Q_m\right) 
\mbox{ is connected}  
\right\},
\end{eqnarray*}
where $|I|$ is the number of the elements of $I.$ 
Then, $\mathcal{Q}^\ast \subset \mathcal{D}^\sharp$
is also a countable set 
and it is dense in 
$\mathcal{D}^\sharp.$
Therefore, $\mathcal{D}^\sharp$ is separable. 
\hspace*{\fill} 
\QED
\medskip

\noindent
\textit{Proof of Theorem \ref{theorem.sect.uniqueness.appl.1}.}
According to Lemma \ref{lemma.sect.uniqueness.1},
the family of open subsets of $\mathcal{D}$ is
countably generated. 
Therefore, 
$\mathcal{B}_{\mathcal{D}}$ is countably generated, 
and $\mathcal{A}_{\mathcal{X}}$ is also countably generated. 
Therefore, the result follows from 
Theorem \ref{theorem.sect.uniqueness.1}. 
\hspace*{\fill} 
\QED

\section*{Acknowledgements}
The author is most grateful to 
Professor Masaaki Tsuchiya for his valuable advice.